\documentstyle[11pt]{siampdg}
\newcounter{machine}
\ifnum\themachine=1 
   \input psfig
\fi

\def\NAT{ I \kern-0.15em N }
\def\COMPLEX{ I \kern-0.50em C }
\def\nat{ I \kern-0.15em N }
\def\complex{ I \kern-0.50em C }

\def\and{ {\rm and} }
\def\half{ {1 \over 2} }

\def\vek#1{{\bf #1}}
\def\vruimte#1#2{{\vrule height #1 depth #2 width 0 pt}}
\def\frac#1#2{\displaystyle{{#1}\over{#2}}}

\def\Im{{\rm Im}}
\def\nopar{\par\noindent}
\def\square#1#2 %maak een vierkant met dimensies #1 en inhoud #2
  {\vbox
    {\hrule
     \hbox
       {\vrule
        \vbox spread#1
          {\vfil
           \hbox spread#1{\hfil#2\hfil}
           \vfil
          }%
        \vrule
        }
     \hrule
     }
  }
\def\sq{\hfill\square{6pt}{} }

\def\beq#1
  {\begin{equation}\label{#1}%\naam{#1}
  }
\def\eeq{\end{equation}}
\def\R{I\!\!R}
\def\br{~|~}
\title{An Introduction to Total Least Squares}
\author{\protect{\parbox[t]{5in}{
  \begin{center}
    P. de Groen\\
    \sl Vrije Universiteit Brussel, \\
    Department of Mathematics,\\
    Pleinlaan 2, B--1050, Belgium \\
    E-mail:~pdegroen@vub.ac.be\\
{\rm This paper has been published in}: \\ 
Nieuw Archief voor Wiskunde, Vierde serie, deel {\bf 14}, 1996, pp. 237-253.
 \end{center}
  }}}
\date{}
\begin{document}

\maketitle

\section{Introduction\label{par1}}
\setcounter{equation}{0}
This (tutorial) paper grew out of the need to 
motivate the usual formulation of a
``Total Least Squares problem'' and to explain the way it is solved
using the ``Singular Value Decomposition''. Although it is an important
generalization of (ordinary) least squares and not 
more difficult to understand,
it is hardly treated in numerical textbooks up to now. 
In the well-known book of Golub \& Van Loan \cite{gvl} 
and in \cite{vanhuffel},
the problem is formulated as follows:
\nopar
\beq{problem}
\matrix{
\mbox{\sl Given a matrix $A\in\R^{m\times n} $ with 
$m>n$ and a vector $\vek b\in\R^m$,}\cr
\mbox{\sl find residuals $E\in\R^{m\times n} $ and $\vek r\in\R^m$
that minimize }\cr
\mbox{\sl the Frobenius norm $\|(\,E\,|\,\vek r\,)\|_F$ 
subject to the condition
$\vek b+\vek r\in Im(A+E)$. }}
\eeq
\nopar
It is proposed as a more natural way to approximate the 
data if both $A$ and $b$
are contaminated by ``errors''.
In our opinion, it is not made clear sufficiently well, 
why this indeed is a natural
generalization of the standard least squares problem
and why it makes sense to study it. On the other hand, 
the classroom note of
Y. Nievergelt \cite{niever} gives a very nice introduction,
but it tells only half of the story in that it considers
(multiple) regression only. 
\par
In this note, we shall give a unified view of
ordinary and total least squares problems and their solution.
As the geometry underlying the problem setting
greatly contributes to the understanding  of the solution, 
we shall introduce
least squares problems and their generalization via 
interpretations in both column space and 
(the dual) row space and we shall use both approaches 
to clarify the solution.
After a study of the least squares approximation for 
simple regression in section \ref{par2},
we introduce the notion of approximation in
the sense of ``Total Least Squares (TLS)'' for 
this problem in section \ref{par3}.
In the next section we consider ordinary and total least 
squares approximations for multiple
regression problems and in section \ref{par5} we 
study the solution of a general overdetermined system
of equations in TLS-sense. In a final section we consider 
generalizations with
multiple right-hand sides and with ``frozen'' columns. We remark
that a TLS-approximation needs not exist in general; however, 
the line (or hyperplane)
of best approximation in TLS-sense for a regression 
problem does exist always.
\par
As numerical algorithms such as the QR-factorization and the
Singular Value Decomposition (SVD) are relatively 
well-known and nicely implemented
in a package like MATLAB, we shall not consider 
numerical algorithms to compute 
the solutions effectively.

\section{Primal vs.~dual approach}
\setcounter{equation}{0}
To make clear how both column- and row-space arguments 
can be used to derive
the solution of a least squares problem, we consider 
least squares in one dimension:
\nopar
\centerline{\sl
Given $m$ points $\{x_i\br i=1,\,\cdots\,,\,m\}$, 
find $z\in\R$ that minimizes the quadratic
functional}
\beq{intro1}
f(z):=\sum_{i=1}^m\,(x_i-z)^2 \,.
\eeq
The function $z\mapsto f(z)$ is a parabola. 
When we shift its center to the average 
$\overline x:={1\over m}\sum_{i=1}^m\,x_i$\,,
\beq{intro1a}
f(z)=\sum_{i=1}^m\,(x_i-z)^2 =\sum_{i=1}^m\,\{\,(x_i-\overline x)^2
 +2(x_i-\overline x)(\overline x - z)+(\overline x-z)^2\,\}\,,
\eeq
we see that the sum of double products vanishes. Hence,
the average $\overline x$
is the unique minimizer.

\nopar
In the dual approach we consider the data as one point
in $\vek x\in\R^m$. The functional $f(z)$ then measures 
the square of the Euclidean distance
to the point $z\vek e$,
\beq{intro2}
f(z)=\|\,\vek x - z\vek e\,\|_2^2\,,~~~\mbox{where}~~~
\vek x :=\left(\matrix{~x_1~\cr x_2 \cr \vdots \cr x_m}\right)
~~\mbox{and}~~
\vek e :=\left(\matrix{~1~\cr 1 \cr \vdots \cr 1}\right)\,.
\eeq

\hrule

\begin{figure}[htb]
\begin{center}
\begin{picture}(280,115)
  \setlength{\unitlength}{.1mm}
  \put(0,0){\line(2,1){800}}
  \thicklines
  \put(200,100){\vector(1,0){500}}
  \put(202,101){\vector(2,1){400}}
  \put(600,300){\vector(1,-2){100}}
  \put(200,60){$\scriptstyle O$}
  \put(700,60){$\scriptstyle \vek x$}
  \put(810,400){$\scriptstyle  span\{\vek e\}$}
  \put(560,320){$\scriptstyle z\,\vek e$}
  \put(650,200){$\scriptstyle \vek x -z\,\vek e$}
\end{picture}
\end{center}
\caption{Vector $\vek x$, its orthogonal projection on $span\{\vek e\}$
     and the residual vector $\vek x -z\,\vek e$ in 
     the dual approach.\label{fig1a}}
\end{figure}
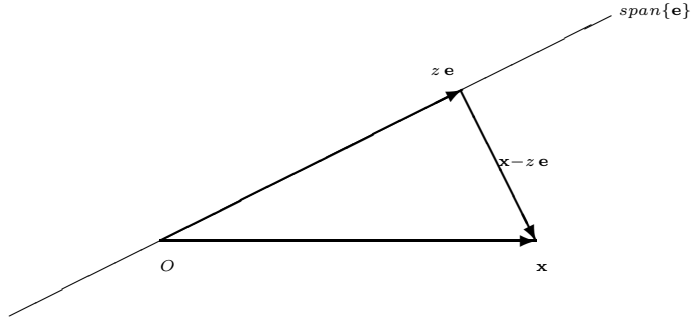

\hrule

\nopar
From fig. \ref{fig1a}, which shows
the plane in $\R^m$ spanned by $\vek x$ and $\vek e$,
we find the orthogonal projection of $\vek x$ on 
$span\{\vek e\}$ as minimizer,
\beq{intro3}
\overline x={\vek x^T\,\vek e \over \vek e^T\,\vek e }=
{1\over m}\sum_{i=1}^m\,x_i\,.
\eeq
We see that both the primal and the dual approach provide 
the solution in different ways. In the primal 
approach we use the fact that linear
terms vanish by a shift
towards the average. In the dual approach we use an orthogonality argument.

\section{Simple regression\label{par2}} 
\setcounter{equation}{0} 
In the plane $\R^2$ we are given $m$ data points (abscissae and ordinates)
\beq{sregres0}
\{(x_i\,,\,y_i)\in\R^2\br i=1,\,\cdots\,,\,m\}
\eeq
that should satisfy
the linear (affine) relation $y(x)=a+bx$; find the parameters $a$ and $b$
that provide a ``best fit'', minimizing the sum of squares of the residuals
\beq{sregres1}
f(a\,,\,b):=\sum_{i=1}^m\,(y_i-a-b\,x_i)^2 \,.
\eeq
We can interpret this as searching the line 
$\ell:=\{(x,y)\in\R^2\br y=a +b\,x\}$  ``nearest'' to the datapoints, 
minimizing vertical distances
and making the tacit assumption that
{\sl model errors in the 
data-model $y=a+bx$ are confined to 
the observed $y$-coordinates}, as depicted in fig. \ref{fig1}.

\ifnum\themachine=0
\begin{figure}[htb]
\hbox{\hskip 37mm \vbox{\vskip 165pt{\special{illustration
tlsfig1a.eps scaled 500}}\vskip -25pt}}
\caption{\label{fig1}
Simple linear regression; {\rm distances are measured along the $y$-axis.}}
%\vskip-15pt
\end{figure}
\else
\begin{figure}[htb]
\vskip -20pt
\hskip43mm\psfig{figure=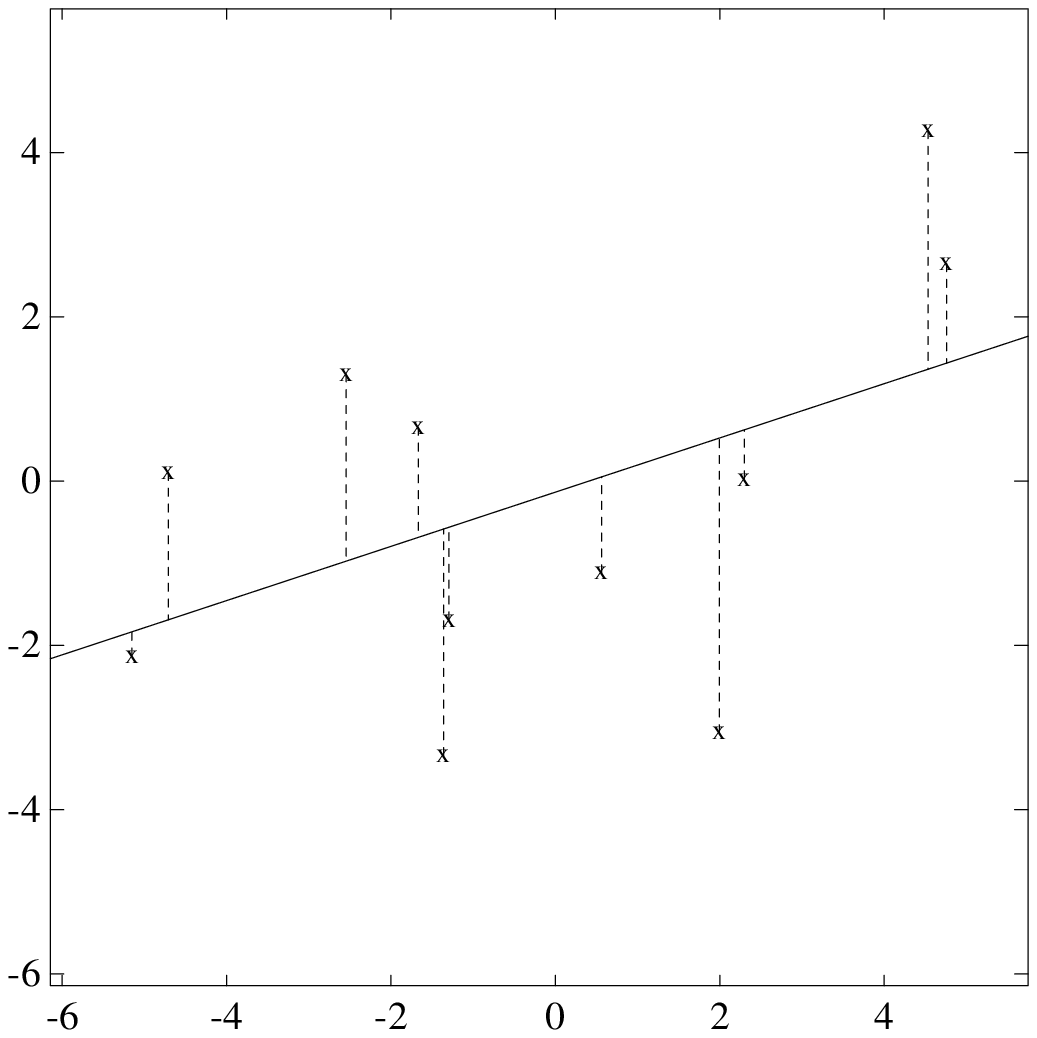,width=210pt,height=160pt}
\vskip-15pt
\caption{\label{fig1}
Simple linear regression; {\rm distances are measured along the $y$-axis.}}
\end{figure}
\fi

\nopar 
Analogously to (\ref{intro1a}) using the centroid
$\overline{ \vek z}:=(\overline x\,,\,\overline y)^T=(
\,{1\over m}\sum_{i=1}^m\,x_i
\,,\,{1\over m}\sum_{i=1}^m\,y_i\,)^T
$  
we rewrite $f$ and find as before, that the double products vanish,
\beq{sregres2}
\begin{array}{r c l}
\displaystyle f(a\,,\,b):=\sum_{i=1}^m\,(y_i-a-b\,x_i)^2 &=&\displaystyle 
\sum_{i=1}^m\,\Big(y_i-\overline y +b\,(x_i-\overline x){\Big )}^2 
+m\,(\overline y - a-b\,\overline x)^2 
\cr ~&\ge&\displaystyle 
\sum_{i=1}^m\,\Big(y_i-\overline y +
b\,(x_i-\overline x)\Big)^2\,,~~~~\forall ~a,\,b\,, \cr
\end{array}
\eeq
with equality if $\overline y = a+b\,\overline x$.
This implies that the centroid is 
located on the line: $\overline{\vek z}\in\ell$.
Eliminating $a$ it remains to minimize a function of $b$ alone, 
which is a parabola. Hence the minimizer of (\ref{sregres1}) is
\beq{sregres3}
b={\sum_{i=1}^m\,(\overline x-x_i)(\overline y-y_i)\over
\sum_{i=1}^m\,(\overline x-x_i)^2}
~~~~\mbox{and}~~~~a=\overline y -b\,\overline x\,.
\eeq

\nopar
{\bf In the dual} approach in $\R^m$ we interpret $x_i$ and 
$y_i$ as components
of vectors $\vek x$ and $\vek y\in \R^m$\,,
\beq{sregres3a}
\vek x :=\left(\matrix{~x_1~\cr x_2 \cr \vdots \cr x_m}\right)
~~~~\vek y :=\left(\matrix{~y_1~\cr y_2 \cr \vdots \cr y_m}\right)
~~~~
\vek e :=\left(\matrix{~1~\cr 1 \cr \vdots \cr 1}\right)~~~~\mbox{and}~~~~
A:=\left(~\vek e\br\vek x\,\right)\in\R^{m\times 2}\,.
\eeq
In this setting the functional $f$ measures the 
square of the distance from $\vek y$ to
a linear combination of $\vek e$ and $\vek x$,
\beq{sregres4}
f(a,\,b)=\|\,\vek y - a\,\vek e-b\,\vek x\,\|_2^2=
\|\,\vek y - A\, {a\choose b} \, \|_2^2\,.
\eeq
As in (\ref{intro3}) it is minimized by the orthogonal 
projection of $\vek y$ on the
span of $\vek x$ and $\vek e$
\beq{sregres5}
f~\mbox{minimal}~~~~\iff ~~~~\vek y - A\, 
{a\choose b} ~\perp~ \mbox{Im}(A)\,.
\eeq
If the rank of $A$ is maximal, the solution
can be computed, see \cite{gvl}, from the {\sl Normal Equations} 
or better by an {\sl Orthogonal Factorization} 
\beq{sregres6}
A^TA {a\choose b}=A^T \vek y~~~~{\rm or~better}~~~~
A=QR~~~\and~~~ R {a\choose b}=Q^T\vek y\,.
\eeq
Otherwise we can use the {\sl Singular Value Decomposition}
\beq{sregres8}
A=U\,\Sigma\,V^T~~~~\and~~~~  {a\choose b}=V\,\Sigma^\dagger\,U^T\vek y\,.
\eeq

\section{Total Least Squares for simple regression\label{par3}}
\setcounter{equation}{0} 
In (\ref{sregres1}) and fig. \ref{fig1} 
we considered the problem of locating
a line nearest to a collection of points, 
where the distance is measured along the 
$y$-axis. It looks ``more natural'' to use 
the (shorter) true Euclidean distance instead,
as drawn in fig. \ref{fig2}, which yields the 
line of {\it Total Least Squares}.

\ifnum\themachine=0
\begin{figure}[htb]
\hbox{\hskip 37mm \vbox{\vskip 160pt{\special{illustration
tlsfig1c.eps scaled 500}}\vskip -25pt}}
\caption{\label{fig2}Line of {\it Total Least Squares}: 
{\rm Model errors
are distributed over 
the $x$- and $y$-coordinates.}}
%\vskip-25pt
\end{figure}
\else
\begin{figure}[htb]
\vskip -20pt
\hskip43mm\psfig{figure=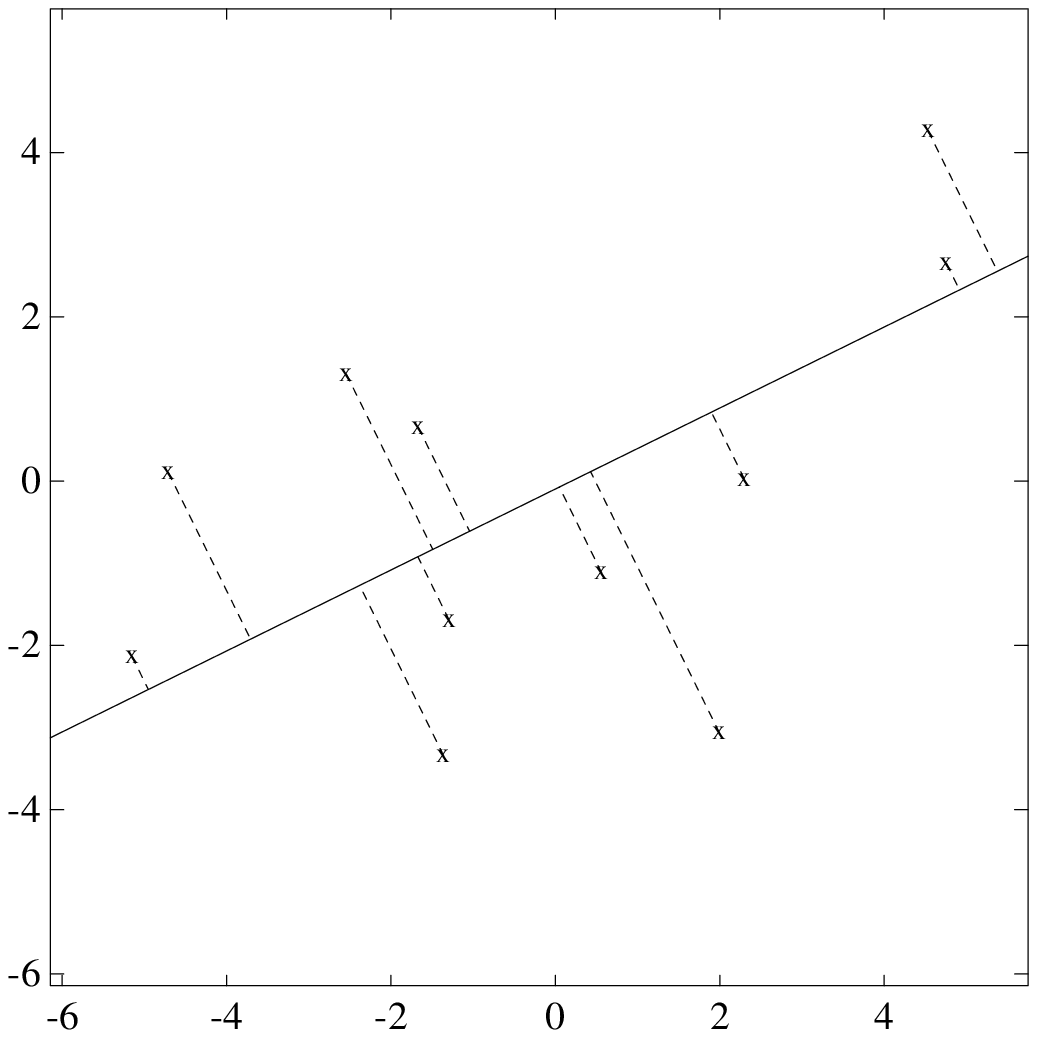,width=210pt,height=160pt}
\vskip-15pt
\caption{\label{fig2}Line of {\it Total Least Squares}: 
{\rm Model errors are distributed over 
the $x$- and $y$-coordinates.}}
\end{figure}
\fi

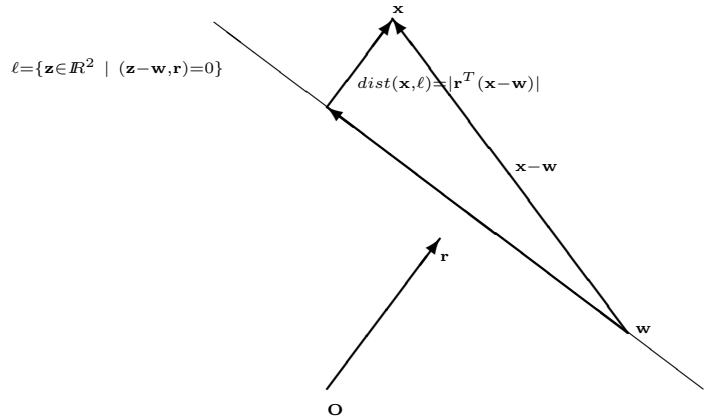
\begin{figure}[htb]
\hrule
\vspace*{15pt}
\begin{center}\setlength{\unitlength}{1pt}
\begin{picture}(280,140)(-130,0)
  \setlength{\unitlength}{.1mm}
  \put(500,0){\line(-4,3){650}}
  \thicklines
  \put(400,75){\vector(-3,4){313}}
  \put(0,0){\vector(3,4){150}}
  \put(0,375){\vector(3,4){87}}
  \put(400,75){\vector(-4,3){400}}
  \put(0,-35){$\scriptstyle\bf O$}
  \put(150,170){$\scriptstyle\vek r$}
  \put(410,75){$\scriptstyle\vek w$}
  \put(-420,420){$\scriptstyle\ell=\{\vek z 
                       \in \R^2~\vert~(\vek z - \vek w,\vek r)=0\}$}
  \put(87,500){$\scriptstyle\vek x$}
  \put(250,290){$\scriptstyle\vek x-\vek w$}
  \put(40,400){$\scriptstyle dist(\vek x,\ell)=
                     {|\vek r^T(\vek x -\vek w)|}$}
\end{picture}
\end{center}
\caption{\label{fig3}The line $\ell$ in the plane is given as the line 
through the vector $\vek w$ orthogonal
to the vector $\vek r$ of unit length. For a given vector $\vek x$ the 
difference vector $\vek x -\vek w$ is drawn
together with its projection along the line $\ell$ and 
its orthogonal complement.}
\vspace*{10pt}
\hrule
\end{figure}

\nopar
So we consider the {\sl Total Least Squares} problem of finding the line 
$\ell$ that minimizes the 
sum of squares of {\bf true} distances:
\beq{tlstwo1}
f(\ell)~:=~\sum_{i=1}^m\,dist(\,(x_i\,,\,y_i)\,,\, \ell\,)^2
\eeq
Instead of asking for a line $y=ax+b$, we use the more symmetric form
\beq{tlstwo1a}
\ell=\{(x,y)\in\R^2\br a+r_1x+r_2y=0\}=\vek w+\vek r^\perp,~~~{\rm with}~~~
\|\vek r\|^2=r_1^2+r_2^2=1,
\eeq 
where $\vek w$ is an arbitrary point on the 
line $\ell$, i.e. $a+r_1w_1+r_2w_2=0$.
With this parametrization of $\ell$ we accept 
the possibility, that $r_2$ may become zero, and hence, that
the line cannot be recast in the form $y=\alpha +\beta x$.
In the description $\ell=\vek w+\vek r^\perp$, where 
$\vek r$ is of unit length,
the distance from a point $\vek z$ to 
$\ell$ is given by, see fig. \ref{fig3},
\beq{tlstwo2}
~~~~~dist(\vek z,\ell)={|\vek r^T(\vek z -\vek w)|}
~~~{\rm where}~~~
\ell=\vek w +\vek r^\perp=\{\vek z \in \R^2 \br 
\vek r^T(\vek z - \vek w)=0\,\}
~~\and~~\|\vek r\|=1\,.
\eeq

\nopar
Hence the TLS problem is to find $\vek r$ and 
$\vek w$ that minimize the functional
\beq{tlstwo3}
I(\vek r,\vek w):=\sum_{i=1}^m~ \left(\vek r^T(\vek z_i -\vek w)\right)^2=
                  \sum_{i=1}^m ~\left(r_1\,(x_i-w_1)+r_2\,(y_i-w_2)\right)^2
\eeq
where
$$
\vek z_i={x_i\choose y_i}~~~~\mbox{and}~~~~
   \vek r={r_1\choose r_2}\,,~~~~\|\vek r\|^2=r_1^2+r_2^2=1\,. 
$$
Making the shift to the centroid, as in (\ref{sregres2}) and (\ref{intro1a}),
we find again, that the sum of double products vanishes, 
\beq{tlstwo4}
\begin{array}{l c l}
I(\vek r,\vek w)&=&\displaystyle
\sum_{i=1}^m ~\left(\vek r^T (\vek z_i-\vek w)\vruimte{1.1em}{.6em}\right)^2 
     \vruimte{0pt}{12pt}
\\
&=&\displaystyle
\sum_{i=1}^m ~\left(\vek r^T (\vek z_i-\overline{\vek z})\right)^2\, 
+\sum_{i=1}^m ~2\,\vek r^T (\vek z_i-\overline{\vek z})\,
                      \vek r^T (\overline{\vek z}-\vek w)
+\,m(\vek r^T (\overline{\vek z}-\vek w))^2\vruimte{1.9em}{1.6em}
\\
&=&\displaystyle
 I(\vek r,\overline{\vek z})
  +\,m(\vek r^T (\overline{\vek z}-\vek w))^2~\ge~ 
 I(\vek r,\overline{\vek z})
\,.\\
 \end{array}
\eeq
Clearly, the centroid $\vek{\overline z}:=(\overline x,\overline y)^T$ 
minimizes the
functional $\vek w\mapsto I(\vek r,\,\vek w\,)$ for every $\vek r\in\R^2$.
This implies, that the minimizing line $\ell=\vek{\overline z}+\vek r^\perp$
passes through the centroid
(as did the line of simple regression) and that we are left with the
reduced minimization problem: 
\\{\sl Find the vector $\vek r$ with $\|\vek r\|_2=1$ minimizing}
\beq{tlstwo5}
I(\vek r,\,\vek{\overline z})=
\sum_{i=1}^m ~\left(r_1\,(x_i-\overline x)+r_2\,(y_i-\overline y)
     \vruimte{1.1em}{.2em}\right)^2=
\|B\vek r\|_2^2=\vek r^T\,B^T\,B\,\vek r\,,
\eeq
where $B\in \R^{m\times 2}$ is the matrix
\beq{tlstwo5a}
B:=\left(\,\vek x-\overline x\,\vek e\br\vek y-\overline y\,\vek e\,\right)=
\left(\matrix{x_1-\overline x & y_{1}-\overline y\cr 
      x_2-\overline x & y_{2}-\overline y\cr \vdots & \vdots\cr
      x_m-\overline x & y_m-\overline y}\right) \,.
\eeq
The problem of minimizing
$\|\,B\,\vek r~\|_2^2$ subject to $\|\,\vek r\,\|_2=1$ is solved by the
Singular Value Decomposition of $B$,
$$
B=U\,\Sigma\,V^T~~~~\mbox{with}~~~~\Sigma=
\left(\matrix{~\sigma_1~&~0~\cr~0~&~\sigma_2~}\right)~~~
\mbox{and}~~~\sigma_1\ge\sigma_2\,.
$$
The solution vector $\vek r$ of (\ref{tlstwo5})
is the right singular vector of $B$ corresponding to the smaller
singular value of $B$\,.
So we conclude:
\begin{list}{\alph{enumi}.}{\leftmargin15pt \usecounter{enumi}}
\item The solution always exists and is given by 
the line through the centroid
orthogonal to the subdominant singular vector of $B$. 
\item
As $r_2$ can be zero, the solution
needs not be expressible in the form $y=\alpha +\beta x$.
\item The solution is unique iff $\sigma_1\ne\sigma_2\,.$
\item The shift (\ref{tlstwo4}) to the centroid $\vek z\in\ell$ is 
the key in finding the solution, as shown in \cite{niever}.
\end{list} 

\ifnum\themachine=0
\begin{figure}[htb]
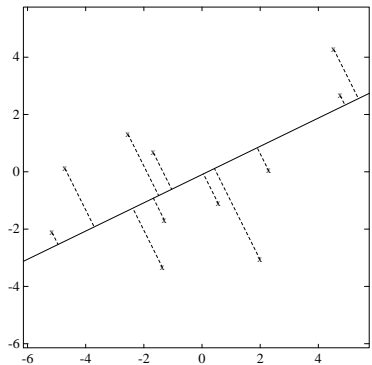

\hbox{\hskip 37mm \vbox{\vskip 160pt{\special{illustration
tlsfig1c.eps scaled 500}}\vskip -15pt}}
\caption{\label{fig2a}Components $(f_i,g_i)$ are the best 
approximations of $(x_i,y_i)$
on the line $a+r_1x+r_2y=0$\,.}
\vskip-10pt
\end{figure}
\else
\begin{figure}[htb]
\vskip -20pt
\hskip43mm\psfig{figure=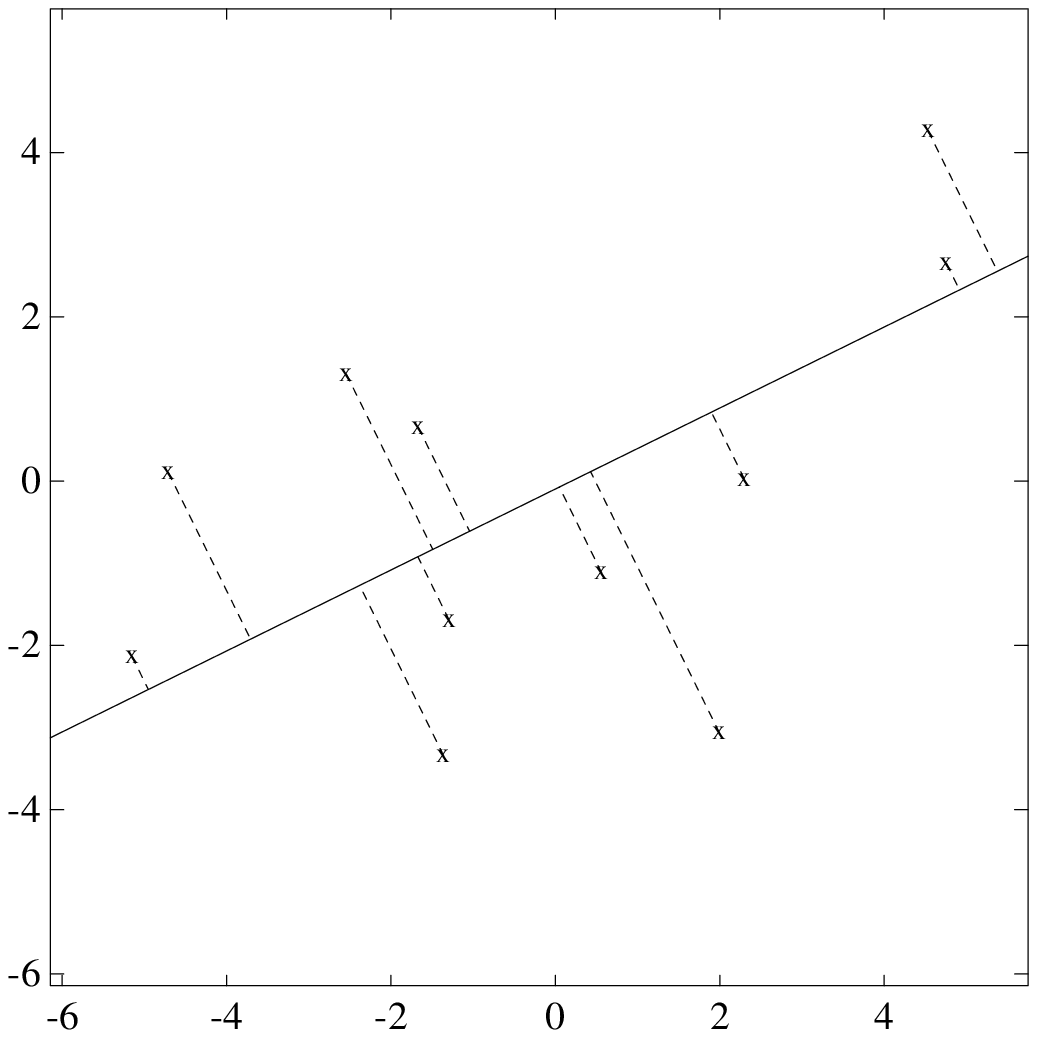,width=210pt,height=160pt}
\vskip-15pt
\caption{\label{fig2a}Components $(f_i,g_i)$ are the best 
approximations of $(x_i,y_i)$
on the line $a+r_1x+r_2y=0$\,.}
\end{figure}
\fi

\nopar
{\bf In the dual formulation} we consider the vectors 
$\vek x$, $\vek y$ and $\vek e$ as in (\ref{sregres3a}) and
we describe the line $\ell$ as in (\ref{tlstwo1a}) by 
$\ell:=\{(\xi,\eta)\br
a +r_1 \xi+r_2 \eta=0\}$.
For $i=1\,\cdots\,m$ we denote by $(f_i,g_i)$ the point on 
$\ell$ nearest to $(x_i,y_i)$, see fig.~\ref{fig2a}, 
and by $(\overline f,\overline g):={1\over m}\sum_{i=1}^m (x_i,y_i)$ we denote
their average.
We define the vectors of first and second components 
$\vek f$, 
$\vek g\in \R^m$,
$$ \vek f:=(f_1\,,\,f_2\,,\,\cdots\,,\,f_m)^T~~~\and~~~
 \vek g:=(g_1\,,\,g_2\,,\,\cdots\,,\,g_m)^T.
$$
These vectors clearly satisfy the relation 
$a\,\vek e+r_1\,\vek f+r_2\, \vek g=0$.
So we can rephrase the minimization problem (\ref{tlstwo1}) as the quest for
vectors $\vek f$ and $\vek g$ that minimize the sum of squares of distances
\beq{tlstwo6}
\begin{array}{r c l}
I(a,\vek r):=\sum_{i=1}^m (x_i-f_i)^2 &+&\sum_{i=1}^m (y_i-g_i)^2~=~
\|\,\vek x-\vek f\,\|^2_2+\|\,\vek y-\vek g\,\|^2_2\cr
& &\mbox{subject to}~~~~
a\,\vek e+r_1\,\vek f+r_2 \,\vek g=0\,,~~~r_1^2+r_2^2=1.\vruimte{1.9em}{0em}
\end{array}
\eeq
Decomposing the vectors in their components in $span\{\vek e\}$ and in the orthogonal
complement $\vek e^\perp$ we obtain
\beq{tlstwo6a}
I(a,\vek r)=\|\,\vek x-\vek f-(\overline x-\overline f)\vek e\,\|^2_2+
\|\,\vek y-\vek g-(\overline y-\overline g)\vek e\,\|^2_2+
m(\overline x-\overline f)^2+m(\overline y-\overline g)^2\,.
\eeq
The contributions from the parts in $span\{\vek e\}$ are minimized by the choice
$\overline f=\overline x$ and $\overline g=\overline y$ and the subsidiary 
condition implies $a+r_1\overline x +r_2\overline y=0$ for that choice. Choosing
$\widetilde \vek f:=\vek f-\overline x\,\vek e$ and $\widetilde\vek g
:=\vek g-\overline y\,\vek e$
we are left with the problem to minimize in $\vek e^\perp$ the functional:
\beq{tlstwo6b}
\|\,\vek x-\overline x\,\vek e - \widetilde\vek f\,\|^2_2+
\|\,\vek y-\overline y\,\vek e -\widetilde\vek g\,\|^2_2~~~~\mbox{subject to}~~~~
r_1\,\widetilde\vek f+r_2\,\widetilde\vek g=\vek 0\,.
%\and~~~ \,.
\eeq
It is not necessary to impose the condition 
$\widetilde \vek f\,,\,\widetilde \vek g\in\vek e^\perp$, since it is automatically
satisfied by the minimizer, because $\vek x-\overline x\,\vek e$
and $\vek x-\overline x\,\vek e$ satisfy this condition.
In matrix notation with 
$B:=\left(\,\vek x -\overline x\,\vek e\br\vek y -
      \overline y\,\vek e\,\right)$
and $E:=\left(\,\widetilde \vek f\br\widetilde \vek g\,\right)$
this minimization problem takes the form
\beq{tlstwo7}
\mbox{minimize}~~~~~\|\,B-E\,\|^2_F~~~~~~\mbox{subject to}~~~~rank(E)=1\,.
\eeq
From the Singular Value Decomposition of $B$, 
$$B=\sigma_1\,\vek u_1\,\vek v_1^T+\sigma_2\,\vek u_2\,\vek v_2^T
~~~~\mbox{we find} ~~~~E=\sigma_1\,\vek u_1\,
    \vek v_1^T\,,~~~~\mbox{provided}~~~~\sigma_1>\sigma_2\,.
$$
Hence the total least squares solution is (as before) given by,
$$
E\,\vek v_2=\vek 0~~~~\mbox{implying}~~~~\vek r =\vek v_2\,.
$$
There is a difference in flavour between both approaches.
Whereas the primal formulation (\ref{tlstwo5}) directly 
produces the minimizing vector,
the dual approach (\ref{tlstwo7})
takes a roundabout. The latter provides a minimizing {\sl matrix} $E$;
the parameters of the line are found only
afterwards as the coefficients in the linear 
combination of the columns of $E$
that equals zero.

\section{Multiple regression\label{par4}}
\setcounter{equation}{0} 
The extension of ordinary and total least squares to multiple regression
is almost straightforward. As most ideas in 2D-regression easily carry
over, we can be brief about it. We are given
the cloud of $m$ datapoints in $\R^n$ (each point consisting
of an ``abscissa'' in $\R^{n-1}$ and an ordinate in $\R$),
\beq{mregres0a}
\{\vek z_i:=(x_1^{(i)},\cdots,x_{n-1}^{(i)},y_i)^T\,\in\,\R^n 
 \br i=1,\,\cdots\,,\,m\}\,,
\eeq
that should satisfy the linear (affine) model $y(x_1\,\cdots\,x_{n-1})=
c_0 + c_1 x_1 + c_2 x_2+\cdots+c_{n-1} x_{n-1}\,.$
In ordinary least squares the parameters are determined by 
minimizing the functional $J$,
\beq{mregres1}
J(\vek c):=\sum_{i=1}^m~(y_i-c_0-c_1x_1^{(i)}-
    \cdots-c_{n-1}x_{n-1}^{(i)})^2\,,
~~~\vek c:=(\,c_{0}\,,\,\cdots\,,\, c_{n-1}\,)^T\,.
\eeq
and we can interpret this as the search for  the best 
fitting hyperplane in $\R^n$\,,
\beq{mregres0}
\{ (x_1\,,\,\cdots\,,\,x_{n-1}\,,\,y)^T\in\R^n\,|
\,y=c_0 + c_1 x_1 + c_2 x_2+\cdots+c_{n-1} x_{n-1}\,\}.
\eeq
As in (\ref{sregres2}), the double products 
vanish by a shift of the center to the centroid,
implying 
$$
J(\vek c)\ge \sum_{i=1}^m~\,\left(y_i-\overline y-
      c_1(x_1^{(i)}-\overline x_{1})-
\cdots-c_{n-1}(x_{n-1}^{(i)}-\overline x_{n-1})\,\right)^2
$$
with equality if 
$\overline y=c_0+c_1\,\overline x_1+\cdots+c_{n-1}\,\overline x_{n-1}$.
Hence, the centroid is in the hyperplane. 
However, more than one unknown parameter is left
and the easy argument of (\ref{sregres3}) cannot be applied directly.
On the other hand, the dual approach (in ``column space'') 
(\ref{sregres3a}-\ref{sregres5}) is 
straightforward and provides the solution
easily. Defining vectors 
$\vek x_k$ and $\vek y\in\R^m$ and the matrix $A\in\R^{m\times n}$,
$$
\vek x_k :=\left(\matrix{~x_k^{(1)}~\cr x_k^{(2)} 
      \cr \vdots \cr x_k^{(m)}}\right),~~~
{\vek y}:=\left(\matrix{y_{1}\cr y_{2}\cr\vdots\cr y_m\cr}\right),~~~
\and~~~
A:=\left(\vek e\,|\,\vek x_1\,|\,\cdots\,|\,\vek x_{n-1}\right)
=\left(\matrix{1 & x_1^{(1)}&\cdots&x_{n-1}^{(1)}\cr 
                 1 & x_1^{(2)}&\cdots&x_{n-1}^{(2)}\cr 
                 \vdots & \vdots&~&\vdots\cr
                 1 & x_1^{(m)}&\cdots&x_{n-1}^{(m)}}\right)
$$
the functional (\ref{mregres1}) takes the form:
\beq{mregres2}
J(\vek c)=\|\,\vek y-c_0\vek e-\cdots-c_{n-1}\vek x_{n-1}\,\|^2=
\|\vek y-A\vek c\|_2^2\,.
\eeq
As in (\ref{intro3}) and (\ref{sregres5}) it is minimized by the 
orthogonal projection of $\vek y$ on the
span of $\vek x_1\,\cdots\,\vek x_{n-1}$ and $\vek e$, i.e. on $Im(A)$,
\beq{mregres3}
f~\mbox{minimal}~~~~\iff ~~~~\vek y - A\, \vek c ~\perp~ \mbox{Im}(A)\,.
\eeq
As before, if the rank of $A$ is maximal, the solution
can be computed from the {\sl Normal Equations} 
or better by an {\sl Orthogonal Factorization}, see \cite{gvl}, 
\beq{mregres4}
A^TA \vek c =A^T \vek y~~~~{\rm or~better}~~~~
A=QR~~~\and~~~ R \vek c=Q^T\vek y\,.
\eeq
Otherwise we can use the {\sl Singular Value Decomposition}
\beq{mregres5}
A=U\,\Sigma\,V^T~~~~\and~~~~  \vek c=V\,\Sigma^\dagger\,U^T\vek y\,.
\eeq

\nopar
{\bf The total least squares approximation} minimizes 
the sum of squares of {\it true} distances.
We do not attribute a special position to the $y$-coordinate and
describe the hyperplane in $\R^n$, as in (\ref{tlstwo1a}), 
by $\vek w + \vek r^\perp$. The functional to
minimize is:
\beq{mregres6}
I(\vek r,\vek w):=\sum_{i=1}^m~ \left(\vek r^T(\vek z_i -\vek w)\right)^2=
                  \sum_{i=1}^m \left(\vek r^T
                       (\vek z_i -\overline{\vek z})\right)^2
             +m(\vek r^T (\overline{\vek z}-\vek w))^2
\eeq
subject to $\|\vek r\|=1\,.$
Since the double products in the second right-hand side
cancel, the centroid (again) is in the hyperplane and it 
minimizes (\ref{mregres6}) for all $\vek r$.
We are left with the reduced minimization problem,
to find $\vek r$ with $\|\vek r\|_2=1$ minimizing
\beq{mregres7}
I(\vek r,\,\vek{\overline z})=
\|B\vek r\|_2^2
\,,~~~~\mbox{with}~~~~
B:=\left(\matrix{x_1^{(1)}-\overline x_1 &\cdots& 
                x_{n-1}^{(1)}-\overline x_{n-1} & y_{1}-\overline y\cr 
      x_1^{(2)}-\overline x_1 &\cdots&
                x_{n-1}^{(2)}-\overline x_{n-1}  & y_{2}-\overline y\cr 
     \vdots &\ &\vdots& \vdots\cr
      x_1^{(m)}-\overline x_1 &\cdots& 
                x_{n-1}^{(m)}-\overline x_{n-1}  & 
                       y_m-\overline y}\right) \,.
\eeq
The solution vector $\vek r$ is the right singular 
vector of $B$ corresponding to the smallest
singular value of $B$. We conclude:
\begin{list}{\alph{enumi}.}{\leftmargin15pt \usecounter{enumi}}
\item A solution always exists; it is given by the hyperplane through 
the centroid and orthogonal to the right singular vector belonging to the
smallest singular value of matrix $B$. It is not expressible in the form 
(\ref{mregres0}) if $r_n=0$.
\item The solution is unique, iff $\sigma_{n-1}>\sigma_n\,.$
\item The shift of (\ref{mregres6}) to the centroid $\vek z\in\ell$ is 
the key in finding the solution.
\end{list} 

\nopar
{\bf In the dual approach} we again consider the hyperplane (\ref{mregres0}),
but now the $y$-coordinate has no special position in the defining equation,
\beq{mregres8a}
\{ (x_1\,,\,\cdots\,,\,x_{n-1}\,,\,y)^T\in\R^n\,|\,
c_0 + c_1 x_1 + c_2 x_2+\cdots+c_{n-1} x_{n-1}+c_n y=0\,\}\,;
\eeq 
instead of $c_n = -1$ we require $\sum_{i=1}^n c_i^2=1$. We choose 
(for each $i$) the point
$(f_1^{(i)}\,,\,\cdots\,,\,f_{n-1}^{(i)}\,,\,g_i)^T$ on this hyperplane
nearest to the datapoint $\vek z_i$, $(i=1\cdots m)$.
The first, second, etc.\ coordinates of these points
form in $\R^m$ the vectors $\vek f_k$ ($k=1\,\cdots\, n-1$) and $\vek g$,
$$ \vek f_k=(f_k^{(1)}\,,\,f_k^{(2)}\,,\,\cdots\,,\,f_k^{(m)})^T~~~\and~~~
 \vek g=(g_1\,,\,g_2\,,\,\cdots\,,\,g_m)^T\,,
$$
which clearly satisfy the relation
$c_0\vek e+c_1\vek f_1+\cdots+c_{n-1}\vek f_{n-1}+c_n\vek g=0$\,.
The minimization of the sum of squares of 
distances from the datapoints $\vek z_i$
to the hyperplane can now be reformulated as the problem of
finding vectors $\vek f_k$ ($k=1\,\cdots\, n-1$) 
and $\vek g$ in $\R^m$ that minimize the functional
\beq{mregres8}
\|\,\vek y-\vek g\,\|^2_2\,+\,\sum_{k=1}^{n-1}\,
  \|\,\vek x_k-\vek f_k\,\|^2_2
~~~~\mbox{ subject to}~~
c_0\,\vek e+c_1\,\vek f_1+\cdots+c_{n-1}\vek f_{n-1}+c_n\,\vek g=0\,,
\eeq
where $\sum_{k=1}^{n}\,c_k^2=1\,.$
As in (\ref{tlstwo6a} -- \ref{tlstwo6b}) we may
restrict this minimization problem to
$\vek e^\perp$ and eliminate the unknown 
$c_0=-c_n\overline y_n-\sum_{k=1}^{n-1}\,c_k\overline x_k$
by orthogonalization w.r.t. $\vek e$; 
essentially this amounts to the same
as the shift to the centroid in the primal approach in $\R^n$. 
So we find the restricted
problem of finding vectors $\vek f_k$ ($k=1\,\cdots\, n-1$) 
and $\vek g$ that minimize
$$
\|\,\vek y-\overline y\,\vek e -\vek g\,\|^2_2\,+
\,\sum_{k=1}^{n-1}\,\|\,\vek x_k-\overline x_k\,\vek e -\vek f_k\,\|^2_2
~~~~\mbox{subject to}~~~~
c_1\,\vek f_1+\cdots+c_{n-1}\vek f_{n-1}+c_n\,\vek g=0\,.
$$
Without imposing it, the minimizing vectors are orthogonal to $\vek e$ 
automatically, as in (\ref{tlstwo6b}). Defining the matrices $B$ and $E$,
$$
B:=\left(\,\vek x_1 -\overline x_1\,\vek e\br\cdots\br
\vek x_{n-1} -\overline x_{n-1}\,\vek e\br
\vek y -\overline y\,\vek e\,\right)~~~
\and ~~~E:=\left(\,\vek f_1\br\cdots\br\vek f_{n-1}\br\vek g\,\right)
$$
we can reformulate the problem as:
\beq{mregres9}
\mbox{minimize}~~~~~\|\,B-E\,\|^2_F~~~~~~\mbox{subject to}~~~~rank(E)=n-1\,.
\eeq
In this form it is easily solved by the SVD.
If $B=\sum_{i=1}^n\,\sigma_i\,\vek u_i\,\vek v_i^T$, then
$E=\sum_{i=1}^{n-1}\,\sigma_i\,\vek u_i\,\vek v_i^T$ 
is a minimizer of (\ref{mregres9}),
which is unique, if $\sigma_{n-1}>\sigma_n\,.$
The coefficients $c_1\,,\,\cdots\,,\,c_n$ determining the hyperplane 
are the coordinates of the right singular vector
$\vek v_n$ as before:
$$
E\,\vek v_n=\vek 0\,,~~~~\Longrightarrow~~~~
\left(\matrix{c_1\cr\vdots\cr c_n}\right)=\vek v_n\,.
$$

\section{General Least Squares\label{par5}}
\setcounter{equation}{0} 
For a given matrix $A\in\R^{m\times n}$ with $m>n$ 
and right-hand side $\vek b\in\R^m$ we
consider the problem to find the  
minimizer $\vek c\in\R^n$ of the functional
\beq{gls1} 
J(\vek c)~:=~\|\,A\,\vek c\,-\,\vek b\,\|_2^2~~~~~
\mbox{with}~~~\vek c:=\left(\matrix{~c_1~\cr\vdots\cr c_n}\right)\,.
\eeq
where
$$A:=\left(\matrix{~a_{1,1}~&~\cdots~&~a_{1,n}~\cr
\vdots& &\vdots\cr a_{m,1}~&~\cdots~&~a_{m,n}}\right)\,\in\R^{m\times n}
~~~~~\and~~~~~
\vek b:=\left(\matrix{~b_1~\cr\vdots\cr b_m}\right)\,\in\R^m~~~(m\ge n)\,, 
$$
The difference with (\ref{mregres2}) is, that $A$ 
needs not contain a column consisting
of all ones.
The solution is obtained by a column space argument as in (\ref{mregres3}),
namely that
$J(\vek x)$ is minimal iff $\vek b - A\, \vek x$ is orthogonal to
$\mbox{Im}(A)$ and it may be computed
by normal equations, QR-factorization or SVD.
\par
What is interesting for the TLS generalization is the interpretation
of (\ref{gls1}) in row space.
We have introduced the TLS approximation in the sections
\ref{par3} and \ref{par4} as the one that minimizes
the sum of squares of the {\em true} distances of 
$m$ points to a hyperplane,
whereas ordinary least squares measures the distances along the $y$-axis.
We can interpret (\ref{gls1}) in this sense.
The rows of the extended matrix $(A\,|\,-\vek b)$ 
define a cloud of $m$ points in $\R^{n+1}\,,$ 
\beq{gls2}
~~~~~\vek z_k :=(\,a_{k,1}\,,\,\cdots\,,\, 
a_{k,n}\,,\,-b_k\,)^T\in\R^{n+1}\,, 
~~~\mbox{such that}~~~
\left(\,\vek z_1\br\cdots\br\vek z_m\,\right)=
\left(\,A\br -\vek b\,\right)^T,
\eeq
to which we try to fit a linear function 
$b(x_1\cdots x_n)=c_1x_1+\cdots+c_nx_n$.
In other words, we look for an $n$-dimensional 
{\em subspace} $\vek{\widehat c}^\perp$
in $\R^{n+1}$
(and not a hyperplane in $\R^n$ as in the regression problem), that
is nearest to the datapoints (\ref{gls2}), minimizing
\beq{gls3}
~~~J(\vek c)=\|\,\left(\,A\br -\vek b\,\right)\,
\left(\matrix{~\vek c~\cr 1}\right)\,\|^2_2~=~
\sum_{k=1}^m~(\,\vek z_k^T\, \vek{\widehat c}\,)^2~~~{\rm where}~~~
\vek{\widehat c}:=\left(\matrix{\vek c\cr 1}\right)=
\left(\matrix{c_1\cr \vdots \cr c_n\cr 1}\right)\in\R^{n+1}\,.
\eeq
In this sum of squares the quantity 
$\vek z_k^T\, \vek{\widehat c}$ measures the distance 
from $\vek z_k$ to $\vek{\widehat c}^\perp$
along the $n\!+\!1$-st coordinate axis. 

\nopar
{\bf The Total Least Squares} approximation for the cloud of
points (\ref{gls2}) minimizes the sum of squares of {\bf true} distances 
to the subspace $\vek{\widehat c}^\perp$.
As the true distance from $\vek z_k$ to the subspace
is given by $\vek z_k^T \vek c /\vek c^T\vek c$\,, 
see (\ref{tlstwo2}), the TLS-approximation
minimizes the functional:
\beq{gls4}
I(\vek c):=
\sum_{k=1}^m~{\left(\,\vek z_k^T\, \vek{\widehat c}\,\right)^2\over 
                     \vek{\widehat c}^T\, \vek{\widehat c}}~=~
{\|\left(\,A\,| -\vek b\,\right)\, \vek{\widehat c}\,\|^2 \over 
                     \vek{\widehat c}^T\, \vek{\widehat c}}
                     ~~~~{\rm where}~~~~
\vek{\widehat c}:=\left(\matrix{\vek c\cr 1}\right)
\eeq
The fuctional $\vek r \mapsto \|(\,A\,| -\vek b\,)\,\vek r\|^2$ 
subject to $\|\vek r\|=1$
is minimal, if $\vek r$ is the right singular vector corresponding to
the smallest singular value of the matrix $(\,A\,| -\vek b\,)$. 
Renormalizing the
last component to $-1$, {\em if possible}, 
provides the solution to the TLS problem
for the overdetermined system of equations $A\vek x=\vek b$.
If the $n\!+\!1$-st component of this right singular
vector is zero, no solution exists to the TLS-problem.
The solution is unique if $\sigma_n>\sigma_{n+1}$.

\nopar
{\bf Interpretation of TLS in Column Space:}
To each point $\vek z_k$ ($k=1\cdots m$) in the cloud (\ref{gls2})
\beq{gls5}
\vek z_k=\left(\matrix{~a_{k,1}~\cr\vdots~\cr a_{k,n}\cr -b_k}~\right)
~~~\mbox{corresponds its best approximation} 
~~~\vek w_k:=
\left(\matrix{~f_{k,1}~\cr\vdots\cr f_{k,n}\cr -g_k}~\right)~
\in~\vek{\widehat c}^\perp\,.
\eeq
The TLS-approximation minimizes the sum of squares of the
distances between the (given) points $\vek z_k$ and the points $\vek w_k$
in the subspace $\vek{\widehat c}^\perp$. We can write
this sum of squares as the Frobenius norm of a matrix, 
if we consider the components $f_{k,j}$
as the elements of a matrix $F\in\R^{m\times n}$, and the components
$g_k$ as the components of a vector $\vek g\in\R^m$.
Hence, TLS minimizes
\beq{gls6}
\sum_{k=1}^m\,\| \vek z_k - \vek w_k\|^2=\|A-F\|^2_F+\|\vek b-\vek g\|^2=
\|(A\,|-b)-(F\,| -g)\|^2_F
\eeq

Since the rows of the matrix $E:=(F\,| -g)\in\R^{m\times (n+1)}$ 
are orthogonal to $\vek{\widehat c}$, the
rank of $E$ is $n$ at most.
In other words,
TLS minimizes
\beq{gls7}
\|\,(\,A\,|-\vek b\,) - E\,\|_F^2
~~~~~\mbox{subject to}~~~~~
E \in \R^{m\times (n+1)}~~~\and~~~rank(E)\le n\,.
\eeq
We may interpret this as the quest for the solution of the
solvable linear system $F\vek c=\vek g$ ``nearest'' to 
the (unsolvable) system $A\vek x=\vek b$, where ``solvable''
means: $\vek g \in {\rm Im}(F)\,$.
\par
The minimization problem (\ref{gls7}) is solved by the SVD.
If $(\,A\,| -\vek b\,)=\sum_{i=1}^{n+1}\,\sigma_i\vek u_i\vek v_i^T\,,$
then $E=\sum_{i=1}^{n}\,\sigma_i\vek u_i\vek v_i^T$ and the required solution
of the TLS-problem is the null-vector $\vek v_{n+1}$ of $E$, 
i.e.\ the right singular vector
$\vek v_{n+1}$ of $(\,A\,|-\vek b\,)$ corresponding 
to the smallest singular value
$\sigma_{n+1}$\,, provided the $n\!+\!1$-st component is non-zero.
As stated at the end of section \ref{par3}, the formulation 
(\ref{gls7}) takes a roundabout
in comparison to the equivalent formulation (\ref{gls4}) 
in that it asks for a minimizing system
of equations, instead of the solution $\vek{\widehat c}$ itself.

\begin{figure}[htb]
\begin{center}
\begin{picture}(300,100)(-150,-55)
  \setlength{\unitlength}{.32mm}
  \put(-100,-50){\line(2,1){200}} % line l
  \put(-100,0){\line(1,0){200}} % horizontal
  \put(0,0){\circle*{3}}
  \put(50,50){\circle*{3}}\put(27,49){$\scriptscriptstyle (1,1)$}
  \put(50,-50){\circle*{3}}\put(22,-53){$\scriptscriptstyle (1,-1)$}
  \put(-50,50){\circle*{3}}\put(-78,49){$\scriptscriptstyle (-1,1)$}
  \put(-50,-50){\circle*{3}}\put(-83,-53){$\scriptscriptstyle (-1,-1)$}
  %\thicklines
  \put(50,50){\line(1,-2){10}}
  \put(50,-50){\line(-1,2){30}}
  \put(50,-50){\line(0,1){100}}
  \put(-25,-7){$\scriptstyle\varphi$}
  \put(53,10){$\scriptstyle\tan\,\varphi$}
  \put(57,40){$\scriptstyle\xi$}
  \put(27,-25){$\scriptstyle\eta$}
\end{picture}
\end{center}\vskip-10pt
\caption{\label{fig4} Example 1: $\xi^2+\eta^2=(1+\tan\varphi)^2\cos^2
           \varphi+(1-\tan\varphi)^2\cos^2\varphi=2$ 
independent on $\varphi$\,.}
%\vskip -30pt
\end{figure}
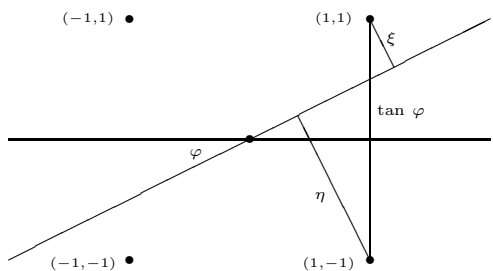

We conclude, that in general a best approximation of the 
overdetermined system
$A\vek x=\vek b$ in TLS-sense may not exist, because we are not
satisfied with the subspace as in a problem of regression; 
we want the equation
for the subspace $b=c_1x_1+\cdots+c_nx_n$ to be explicit w.r.t.\ $b$.
Furthermore, the solution is not necessarily unique.
We shall illustrate this by two examples.
\nopar
{\bf Example 1:} Consider the cloud of 4 points in $\R^2$:
$$
(1,1)\,,~~(-1,1)\,,~~(1,-1)\,,~~\and~~(-1,-1)
$$
The LS-approximation is the horizontal line $\{(x,y)\br y=0\}$.
The TLS-approximation makes the SVD of the matrix $B$,
$$
B:=\left(\matrix{~1~&~1~\cr~1&-1~\cr-1~&~1\cr-1~&-1~}\right)=
   \left(\matrix{~\half ~&~\half ~& \sqrt{\half}&0\cr
                 ~\half &-\half~ & 0& \sqrt{\half}\cr
                  -\half~ &~\half& 0& \sqrt{\half} \cr
                  -\half ~&-\half&\sqrt{\half}&0 ~
              }\right)
   ~\left(\matrix{~2~&~0~\cr0&2\cr0&0\cr0&0}\right)~
\left(\matrix{~1~&~0~\cr0&1}\right)\,.
$$
As both singular values are equal, there is no unicity; 
every line through the origin 
provides a solution, as shown in fig.\ \ref{fig4}. 
The sum of squares of distances
from the points to a line with slope $\tan\phi$ is 
independent of the slope. 
\medskip

\hrule

\hrule

\nopar
{\bf Example 2:} Solve the following problem in LS-sense and TLS-sense:
$$
\left(\matrix{~1~&~0~\cr 0&0\cr0&0}\right)~{x \choose y}=
\left(\matrix{~1~\cr1\cr 1}\right)
$$
The normal equations for the LS-approximation are:
$$
\left(\matrix{~1~&~0~\cr 0&0}\right)~{x \choose y}=
\left(\matrix{~1~\cr 0}\right)~~~~~\Longrightarrow
~~~~x=1~~~\and~~~y~\mbox{undetermined}\,.
$$
The SVD for TLS-problem is:
$$
\def\ha{\scriptscriptstyle{1\over 2}}
\def\mha{\scriptscriptstyle{-{1\over 2}}}
\def\rha{\scriptscriptstyle{1\over\sqrt 2}}
\def\mrha{\scriptscriptstyle{-{1\over\sqrt 2}}}
\def\rpr{\scriptscriptstyle{\sqrt{2+\sqrt 2}}}
\def\rmr{\scriptscriptstyle{\sqrt{2-\sqrt 2}}}
B~=~\left(\matrix{~1~&0&~1~\cr0&0&1\cr0&0&1}\right)~=~
\left(\matrix{ \rha & \rha & 0 \cr \ha & \mha & \mrha \cr \ha 
       & \mha & \rha }\right) 
\left(\matrix{ \rpr & 0 & 0 \cr0&\rmr&0\cr0&0&0}\right) 
\left(\matrix{ \ha\rmr & 0 & \ha\rpr \cr \ha\rpr & 0 & 
       \mha\rmr \cr0&1&0}\right)\,.
$$
The smallest singular value is 0\,.
However, the $3^{\rm rd}$ component of the 
corresponding right singular vector 
$(0\,,\,1\,,\,0)^T$ is 0 as well, such that no TLS-solution exists!

\section{Generalizations: (a) Multiple RHS\label{par7}} 
\setcounter{equation}{0} 
In ordinary
least squares there is no difference between the treatment of one and
multiple right-hand sides (RHS). In Total Least Squares the 
column space of the matrix
is bent towards the RHS. If there are given several RHS's,
we can treat each of them separately and compute the SVD 
of an extended matrix
for each RHS. In a different approach we can try to bend the matrix 
to all RHS's collectively. So we consider the problem:
given $A\in\R^{m\times n}$ ($m\ge n+p$) and 
$B\in\R^{m\times p}$ find $X\in\R^{n\times p}$
that solves the overdetermined system of equations $A\,X=B$ in TLS-sense.
By analogy to (\ref{gls6}) we have to
find the solution $X$ of a solvable matrix equation 
$F\,X=G$ (i.e. $\Im(G)\subset\Im(F)$\,)
nearest to $A\,X=B$; we have to minimize
\beq{gentls1}
~~~\|\,A-F\,\|_F^2~+~\|\,B-G\,\|_F^2
~~~\mbox{subject to}~~~F \in \R^{m\times n}\,,~G 
    \in \R^{m\times p}~\and~F\,X=G\,.
\eeq
Otherwise stated, find an approximation
$E=(\,F\br G\,) \in \R^{m\times(n+p)}$ to $(\, A \br B \,)$, such that
\beq{gentls2}
\|\,(\,A\br B\,)-E\,\|_F^2~~~~\mbox{is minimal subject to }~~~~rank(E)=n\,.
\eeq
The solution of (\ref{gentls2}) is constructed by 
making the SVD of $(\,A\br B\,)$:
\beq{gentls3}
\def\stapelup#1#2{\matrix{{\scriptscriptstyle #2}\cr #1}}
\def\stapel#1#2{\matrix{#1\cr{\scriptscriptstyle #2}}}
(\, A \br B \,)=U\,\Sigma\,V^T=
\left(\, \stapelup{U_1\vruimte{1em}{1em} }{(m\times n)}
        ~\vrule~\stapelup{ U_2\vruimte{1em}{1em} }{(m\times p)}\,\right)~
\left(\matrix{~\stapelup{\Sigma_{1}}{ (n\times n)}&~0~\cr  
               \vruimte{1em}{0em}~&~\cr
              0&\stapel{\Sigma_{2}}{ (p\times p)}}\right)~
\left(\matrix{~\stapelup{V_{1,1}}{(n\times n)}&~\stapelup{V_{1,2}}
              { (n\times p)}~\cr 
               \vruimte{1em}{0em}~&~\cr
               \stapel{V_{2,1}}{ (p\times n)}&\stapel{V_{2,2}}
            { (p\times p)}}\right)^T\,.
\eeq
\nopar
{\bf Theorem.} If we assume:
\begin{list}{\alph{enumi}.}{\leftmargin15pt \usecounter{enumi}}
\item $rank(V_{2,2})=p$\,,
\item $\Sigma=diag(\sigma_1\,,\,\cdots\,,\,\sigma_n\,,\,\sigma_{n+1}
\,,\,\cdots\,,\,\sigma_{n+p}\,)$ with $\sigma_j\ge\sigma_{j+1}$ 
and $\sigma_n\ne\sigma_{n+1}$\,,
\end{list}
then the TLS problem (\ref{gentls2}) has the unique solution 
$X=-V_{1,2}\,V_{2,2}^{-1}\,.$
\medskip\nopar
{\bf Proof}: From (\ref{gentls3}) and the assumption 
$\sigma_n>\sigma_{n+1}$ it follows, that 
the best $rank ~n$ approximation\footnote{see \cite{gvl} 
theorem 2.5.2} of $(A\br B)$ in the Frobenius norm
is given by $E$,
\beq{gentls4}
E:=\left(\, U_1
        \br U_2\,\right)~
\left(\matrix{~\Sigma_{1}~&~0~\cr
              0&0}\right)~
\left(\matrix{~V_{1,1}~&~V_{1,2}~\cr
               V_{2,1}&V_{2,2}}\right)^T~
=U_1\,\Sigma_1\,\left(\, V_{1,1}^T
        \br V_{1,2}^T\,\right)=(F\br G)\,,
\eeq
where $F:=U_1\,\Sigma_1\, V_{1,1}^T$ and $G:=U_1\,\Sigma_1\,V_{1,2}^T\,$.
The orthogonality of the columns of $V$ implies
$$
{V_{1,1}\choose V_{2,1}}^T\,{V_{1,2}\choose V_{2,2}}=(0)
~~~\mbox{and hence}~~~E\,\displaystyle{V_{1,2}\choose V_{2,2}}=
    F\,V_{1,2}+G\, V_{2,2}=(0)\,.
$$
Under the assumption $rank(V_{2,2})=p$ we may conclude, 
that $X:=-V_{1,2}\,V_{2,2}^{-1}$
solves the approximate equation $FX=G$\,. \sq\nopar
{\bf (b) Fixed columns:}
In section \ref{par2} we have introduced the simple (bivariate)
regression problem and we have shown that it is solved in LS-sense
by the LS-solution of the overdetermined system of equations
$A{a\choose b}=(\vek e \br\vek x){a\choose b}=\vek y$ 
(cf. \ref{sregres4}). However, as explained in
section \ref{par5}, the TLS-solution of this overdetermined 
system of equations 
is derived from the SVD of the matrix 
$(\vek e\br\vek x\br\vek y)\in\R^{m\times 3}$.
This differs from the TLS-solution of the regression problem, 
which is derived from the
SVD of 
$B:=(\vek x -\overline x\vek e\br\vek y-\overline y\vek e)\in\R^{m\times 2}$, 
cf. eq.~(\ref{tlstwo5a}).
The reason for this difference is, 
that the formulation of the regression problem as an
overdetermined set of equations
$A{a\choose b}=\vek y$ hast lost its geometric 
interpretation as a line $y=a+bx$
in the $(x,y)$-plane. In the LS-solution this makes no difference
since all uncertainty is put in the $\vek y$-column.
However, TLS for $A{a\choose b}=\vek y$ puts uncertainty in all three columns 
$\vek e$, $\vek x$ and $\vek y$, although in the regression problem
there is no reason to postulate uncertainty in the ``constant term''.
The TLS-solution of the regression problem can be regained from 
$A{a\choose b}=\vek y$
if we ``freeze'' the first column of $A$ and put uncertainty in the columns
$\vek x$ and $\vek y$ only as in eq.~(\ref{tlstwo6}). The solution is
obtained by orthogonalization w.r.t. the frozen column $\vek e$.
\par
This motivates the study of TLS-problem for $A\,X=B$ with 
frozen columns, see \cite{ghs}, where uncertainty 
is postulated in a part of the columns of $A$ 
(LS is a special case, all columns of the matrix being frozen!). 
So we assume that the matrix
$A$ is partitioned in a frozen part $A_1\in\R^{m\times j}$ and
a part  $A_2\in\R^{m\times k}$ containing some uncertainty with $j+k=n$.
Given a right-hand side $B\in\R^{m\times p}$ with $m\ge j+k+p$\,,
we seek matrices $X_1\in\R^{j\times p}$ and $X_2\in\R^{k\times p}$,
such that
\beq{gentls5}
A_1\,X_1+A_2\,X_2~=~B~~~~~~\mbox{in TLS-sense w.r.t.}~~
A_2~~\and~~B~~~\mbox{keeping $A_1$ fixed.}
\eeq
More precise, minimize among all 
$C\in\R^{m\times k}$ and $D\in\R^{m\times p}$
\beq{gentls6}
\|\,A_2-C\,\|_F^2~+~\|\,B-D\,\|_F^2
~~~~~~\mbox{subject to}~~~~~A_1\,X_1+C\,X_2 =D\,.
\eeq
or otherwise said, subject to the condition 
$rank(A_1\br C\br D)=j+k=n$.
\nopar
Guided by the idea of (\ref{tlstwo6}), where we orthogonalized w.r.t.~the
frozen column, we find the \\{\bf solution}:
\beq{gentls7}
\begin{array}{l}
\mbox{a. Orthogonalize columns of $A_2$ and $B$ 
w.r.t. columns of $A_1$}\hspace{5em}\cr
\mbox{b. Solve TLS-problem in the orthogonal complement 
${\rm Im}(A_1)^\perp$\,.}
\end{array}
\eeq

\nopar 
{\bf Proof:} If $A_1$ is of full column rank ($rank(A_1)=j$), 
we make the QR-factorization 
$$A_1=U\,{R_1\choose 0}~~~~ 
\mbox{ with}~~~~ U\in\R^{m\times m}\mbox{ orthogonal}~~~
\and ~~~R_1\in\R^{j\times j}\,.
$$
Because the Frobenius norm is orthogonally invariant, 
the functional (\ref{gentls6})
is equal to 
\beq{gentls8a}
\|\,U^T A_2-U^T C\,\|_F^2~+~\|\,U^T B-U^T D\,\|_F^2\,.
\eeq
Partitioning the matrices in parts consisting of the topmost $j$ rows and 
the remaining $m-j$ rows respectively,
\beq{gentls8b}
{A_{12}\choose A_{22}}:=U^T\,A_2,~~~
{B_1\choose B_2}:=U^T\,B
,~~~{C_1\choose C_2}:=U^T\,C, ~~~{D_1\choose D_2}:=U^T\,D
%~~~\mbox{ with}~~~{\mbox{$j$ rows}\choose\mbox{$m-j$ rows}}
\,,
\eeq
we can rewrite the functional as
\beq{gentls8}
\|\,A_{12}-C_1\,\|_F^2~+~\|\,B_1-D_1\,\|_F^2~+~\|\,A_{22}-C_2\,\|_F^2
~+~\|\,B_2-D_2\,\|_F^2\,.
\eeq
It has to be minimized subject to the equations
$R_1\,X_1+C_1\,X_2=D_1$ and $C_2\,X_2=D_2\,$. If $X_2$ is known, and if
we choose
$A_{12}=C_1$ and $B_1=D_1$, the first two terms in 
(\ref{gentls8}) vanish and 
$X_1$ can be solved from the equation $R_1\,X_1+C_1\,X_2=D_1$.
Hence it suffices to minimize 
\beq{gentls9}
\|\,A_{22}-C_2\,\|_F^2~+~\|\,B_2-D_2\,\|_F^2~~~~
\mbox{subject to}~~~~C_2\,X_2=D_2\,.
\eeq
This is solved as eq.~(\ref{gentls2}) by the SVD of $(C_2\br D_2)$.
\par
If $A$ is not of full column rank ($rank(A_1)=r<j$), we use the SVD of $A_1$:
$$A_1=U\,\left(\matrix{\Sigma_1 &0\cr0&0}\right)\,
{V^T_1\choose V^T_2}~~~~ 
\mbox{ with}~~~~ U\in\R^{m\times m},~~~
\Sigma_1\in\R^{r\times r},~~~V_1\in\R^{j\times r},
~~~V_2\in\R^{j\times(j- r)}\,.
$$
With the same partitioning as in (\ref{gentls8b}), 
but now with the $r$ topmost rows in the
upper parts and the remaining $m-r$ rows in the lower 
parts, we arrive at the minimization of
(\ref{gentls8}) subject to the conditions
\beq{gentls10}
\Sigma_1\,V_1^T\,X_1+C_1\,X_2=D_1~~~~ \and ~~~~ C_2\,X_2=D_2\,.
\eeq
Choosing
$A_{12}=C_1$ and $B_1=D_1$ and solving $X_2$ from (\ref{gentls9}) 
we can solve $V_1^T\,X_1$ from (\ref{gentls10}). 
This makes the first two terms 
in (\ref{gentls8}) zero, such that the problem again is 
reduced to the form (\ref{gentls2}).
As in standard LS-problems in which the matrix 
is not of full column rank,
the part $X_1$ is not uniquely defined; we may add to it 
any linear combination of the columns of $V_2$\,.\sq


\begin{thebibliography}{9}
\bibitem{ghs}
\begin{quote}
     G.H. Golub. A. Hoffman \& G.W.Stewart,
     {\it A generalization of the Eckart-Young-Mirsky matrix approximation
theorem},
 Linear Algebra and Its Applications,
 {\bf 88/89}, pages 317--327, 1987.
\end{quote}

\bibitem{gvl}
\begin{quote}
     G.H. Golub \& C. Van Loan,
     {\it Matrix Computations}, 
     The Johns Hopkins University Press, 2nd ed. 1989.
\end{quote}

\bibitem{niever}
\begin{quote}
     Y. Nievergelt,
     {\it Total least squares: state of the art 
          regression in numerical analysis}, 
     SIAM Review, {\bf 36} pp. 258 - 264, 1994.
\end{quote}

\bibitem{vanhuffel}
\begin{quote}
     S. Van Huffel \& J. Vanderwalle,
     {\it The Total Least Squares Problem: 
           Computational Aspects and Analysis}, 
     SIAM, Philadelphia, PA, 1991.
\end{quote}


\end{thebibliography}
\end{document}